\documentclass[14pt, reqno]{amsart}
\usepackage{graphicx,psfrag}
\usepackage[psamsfonts]{amssymb}
\usepackage{pdfsync}
\usepackage{amscd}
\usepackage{amsmath}
\usepackage{mathrsfs}
\usepackage{extarrows}
\usepackage{stmaryrd}
\usepackage[small,nohug,heads=littlevee]{diagrams}
\usepackage{comment}
\usepackage{enumerate}
\usepackage{xfrac}
\usepackage{cases}
\diagramstyle[labelstyle=\scriptstyle]

\theoremstyle{plain}
    \newtheorem{theorem}{Theorem}[section]
    \newtheorem{lemma}[theorem]{Lemma} 
    \newtheorem{corollary}[theorem]{Corollary}
    \newtheorem{proposition}[theorem]{Proposition}
    
    \theoremstyle{definition}
    \newtheorem{definition}[theorem]{Definition}
    
    \newtheorem{example}[theorem]{Example}

    \newtheorem{remark}[theorem]{Remark}

\theoremstyle{remark}

\numberwithin{equation}{section}

    \newcommand{\R}{\mathbb{R}}
    
    \newcommand{\C}{\mathbb{C}} 
    
    \newcommand{\Z}{\mathbb{Z}}

        \newcommand{\g}{\mathfrak{g}} 
         \newcommand{\h}{\mathfrak{h}} 
         \newcommand{\p}{\mathfrak{p}}

     \newcommand{\pt}{\text{\it pt}}

\title{A $K$-homological Approach  to The Quantization Commutes with Reduction Problem}

\author{Yanli Song}

\date{\today}

\begin{document}

\abovedisplayskip=2pt
\belowdisplayskip=2pt

\maketitle

\begin{abstract}
In \cite{Kasparov88}, Kasparov defined a distinguished $K$-homology fundamental class, so called the Dirac element. We prove a localization formula for the Dirac element in $K$-homology of crossed product of $C^{*}$-algebras. Then we define the quantization of  Hamiltonian $G$-spaces as a push-forward of the Dirac element. With this, we develop a $K$-homological approach to the quantization commutes with reduction theorem.  
\end{abstract}

\tableofcontents

\section{Introduction}
The purpose of this paper is to study the quantization commutes with reduction principle, which was formulated  by Guillemin and Sternberg \cite{Guillemin82}. There are a number of approaches to the solution such as the geometric proof in \cite{Meinrenken98} and analytic proofs  in \cite{Zhang98, Paradan01}. By making compromises to both geometry and analysis, we express the quantization commutes with reduction problem in the framework of $K$-homology.  Besides providing a simplified proof,  the techniques we used can be adapted to more general situations such as the quasi-Hamiltonian spaces \cite{Alekseev98}, which is our original motivation.

Let $G$ be a compact Lie group, and $M$ a compact $G$-manifold with a $G$-equivariant Spin$^{c}$-structure. We denote by $S_{M}$ its spinor bundle. Regarded as a Hilbert bi-module, the bundle $S_{M}$ induces a push-forward map in  $K$-homology
\begin{equation}
\label{push-forward}
\phi_{S_{M}} :KK_{0}^{G}(C_{\tau}(M), \C) \to KK_{0}^{G}(\C,\C),
\end{equation}
where $C_{\tau}(M)$ is the $C^{*}$-algebra of continuous sections of the complex Clifford bundle $\mathrm{Cliff}(TM)$. The $K$-homology on the right-hand side is isomorphic to the representation ring $R(G)$ and the one on the left-hand side contains a distinguished element:
\[
[d_{M}] \in KK_{0}^{G}(C_{\tau}(M), \C),
\] 
the \emph{Dirac element} \cite{Kasparov88}, defined using the de Rham differential operator $d$. The push-forward index of the Spin$^{c}$-manifold $M$ is defined by:
\[
 \phi_{S_{M}}([d_{M}]) \in KK_{0}^{G}(\C,\C) \cong R(G).
\]
It coincides with the $G$-index of the associated Spin$^{c}$-Dirac operator on $M$.

We are particularly interested in the case when $G$ is a compact connected Lie group and $(M, \omega)$ is a compact symplectic $G$-manifold.  We also assume that $(M, \omega)$ is \emph{pre-quantizable}, that is, there exists a $G$-equivariant Hermitian line bundle $L$ with a Hermitian connection $\nabla^{L}$ such that 
\[
\frac{\sqrt{-1}}{2\pi} (\nabla^{L})^{2} = \omega.
\] 
Let $S_{M}$ be the spinor bundle obtained by coupling the pre-quantum line bundle $L$ with the spinor bundle determined by $\omega$. We can formulate the \emph{quantization} of the Hamiltonian $G$-space $(M, \omega)$ as the image of the Dirac element $[d_{M}]$ under the push-forward map $\phi_{S_{M}}$ \cite{MR2989614}. We denote it by $Q(M, \omega) \in R(G)$.

When the action of $G$ on $M$ is Hamiltonian, there is a \emph{moment map}
\[
\mu : M \to \g^{*},
\]
where $\g^{*}$ is the dual of the Lie algebra $\g$ of $G$. If $0$ is a regular value of $\mu$, then $G$ acts on $\mu^{-1}(0)$ locally freely. For simplicity, we assume that the action is free. We can define the symplectic reduction \cite{Weinstein79}:
\[
M_{0} = \mu^{-1}(0)/G, \hspace{5mm} L_{0} = (L|_{\mu^{-1}(0)}) /G,
\]
where $M_{0}$ is a symplectic manifold with a symplectic 2-form $\omega_{0}$ induced from $\omega$. As before, the quantization of the reduced space $(M_{0}, \omega_{0})$ is defined as the push-forward of the Dirac element $[d_{M_{0}}]$ under $\phi_{S_{M_{0}}}$, denoted by $Q(M_{0}, \omega_{0}) \in \Z$. 
 
The quantization commutes with reduction theorem states as follow:
\begin{theorem}
Let $G$ be a compact connected Lie group, and $(M, \omega)$ a compact pre-quantized Hamiltonian $G$-manifold, with moment map $\mu : M \to \g^{*}$. If 0 is a regular value of $\mu$, then
\label{[Q,R]=0}
\[
Q(M, \omega)^{G} = Q(M_{0}, \omega_{0}). 
\]
\end{theorem}

\begin{remark}
Guillemin and Sternberg \cite{Guillemin82} had proved a quantization commutes with reduction theorem for  K\"{a}hler manifolds. They conjectured the above result for Dirac operators on general symplectic manifolds, see \cite{Guillemin94}. For the case that $G$ is abelian, the conjecture was solved by Vergne\cite{Vergne94} and Meinrenken \cite{Meinrenken96} independently. The general case was first proved by Meinrenken  \cite{Meinrenken98} using symplectic surgery, and re-proved by Tian-Zhang \cite{Zhang98} using analytic techniques and by Paradan \cite{Paradan01} using transversally elliptic operators. After these results, the quantization commutes with reduction theorem has been generalized in various directions. For example, Ma-Zhang \cite{Zhang09,Zhang14} solved an extended version of Vergne's conjecture \cite{MR2334206} on compact Lie group acting on non-compact manifolds. Later on, Paradan\cite{MR2869441} gave a different proof.  For the case when both manifolds and Lie groups are non-compact, one can find related works by Hochs, Landsman, Mathai and Zhang \cite{Hochs08,Hochs09,Landsman05, Mathai10}. At last Paradan and Vergne \cite{MR3258255} proved a quantization commutes with reduction formula for Spin$^{c}$-manifold very recently. 
\end{remark}

We now turn to an introduction of our method. We first extend the push-forward map (\ref{push-forward}) to the  $K$-homology of crossed products of $C^{*}$-algebras:
\begin{equation}
\label{p-2}
\hat{\phi}_{S_{M}} :KK_{0}(C^{*}(G, C_{\tau}(M)), \C) \to KK_{0}(C^{*}(G),\C).
\end{equation}
The target of this map can be identified with $\hat{R}(G) = \mathrm{Hom}_{\Z}(R(G), \Z)$, that is,  the group of possibly infinite $\Z$-linear combination $\sum a_{\gamma} \cdot V_{\gamma}$ of irreducible representations.

Following an idea of Witten \cite{Witten82}, we perturb the de Rham differential operator $d$ by a $G$-invariant vector field $X$ and define a perturbed Dirac element. We show  that the corresponding perturbed Dirac element localizes to a $G$-invariant open neighborhood $U$ of the vanishing set $\{ X = 0 \} \subseteq M$.

In the case of Hamiltonian $G$-space $(M, \omega)$, our method is a variation of Paradan's approach \cite{Paradan01}. We take  $X$ to be the Hamiltonian vector field of the norm square of the moment map and  $U$ a $G$-invariant small open neighborhood of $\{X = 0\}$. The submanifold $U$ is disconnected in general. By the symplectic cross-section theorem, every connected component $U_{\alpha}$ has the geometric structure:
\begin{equation}
\label{extra-2}
U_{\alpha} \cong G \times_{H_{\alpha}} V_{\alpha},
\end{equation}
where $H_{\alpha} \subseteq G$ is a sub-group of $G$ of equal rank, and $V_{\alpha}$ is Hamiltonian $H_{\alpha}$-manifold. The tangent bundle $TU_{\alpha}$ splits equivariantly into ``vertical direction" and ``horizontal direction". Hence, so is the Spin$^{c}$-structure on $U_{\alpha}$ and its induced push-forward map.  

To prove the quantization commutes with reduction theorem, we need to show that the push-forward index of a connected component $U_{\alpha}$ doesn't contribute to the trivial $G$-representation unless it contains $\mu^{-1}(0)$. 
By the splitting (\ref{extra-2}) and the Dirac induction formula given by Kostant \cite{Kostant99}, it is enough to show that the push-forward index of the smaller manifold $V_{\alpha}$ has zero multiplicity at the trivial $H_{\alpha}$-representation. 

We then repeat the localization to the Hamiltonian $H_{\alpha}$-manifold $V_{\alpha}$ and reduce the problem to the case when the group is abelian and the manifold is the total space of some vector bundle by induction. With some computation in Section 4, we show that the proof of the quantization commutes with reduction theorem can be localized to a small neighborhood of $\mu^{-1}(0)$.

$\mathbf{Acknowledgments}$: The author benefited a lot from discussions from Victor Guillemin, Nigel Higson and Eckhard Meinrenken. The author also would like to thank Michele Vergne and Paul-Emile Paradan for reading the first version of this paper and giving some useful advice.

\section{Analytic K-Homology of $C^{*}$-algebras}
In this section we give a brief review on Kasparov's approach to analytic $K$-homology \cite{Kasparov80, Kasparov88}. 

\begin{definition}
Let $A$ be a separable $\Z_{2}$-graded C$^{*}$-algebra. A \emph{Fredholm module} over $A$ is a pair $(\mathcal{H}, \mathcal{F})$:
\begin{itemize}
\item
 a separable $\Z_{2}$-graded Hilbert space $\mathcal{H}$ with an even $*$-homomorphism:  
\[
\rho : A \to \mathbb{B}(\mathcal{H});
\]
\item
 an odd operator $\mathcal{F} \in \mathbb{B}(\mathcal{H})$ such that for all $a \in A$
\[
(\mathcal{F} - \mathcal{F}^{*}) \rho(a) \in \mathbb{K}(\mathcal{H}), \ (\mathcal{F}^{2} - 1) \rho(a) \in \mathbb{K}(\mathcal{H}), \  [\rho(a), \mathcal{F}]  \in \mathbb{K}(\mathcal{H}).
\]
Here $\mathbb{K}(\mathcal{H})$  and $\mathbb{B}(\mathcal{H})$ denote the compact operators and the bounded operators on $\mathcal{H}$ respectively. 
\end{itemize}

\end{definition}

\begin{definition}
\label{equ-1}
An \emph{operator homotopy} between two Fredholm modules over $A$: $(\mathcal{H}, \mathcal{F})$ and $(\mathcal{H}, \mathcal{F}^{'})$, is a norm continuous path of Fredholm modules over $A$
\[
t \mapsto (\mathcal{H}, \mathcal{F}_{t}), \hspace{5mm} t \in [0,1],
\]
with $\mathcal{F}_{0} = \mathcal{F}, \mathcal{F}_{1} = \mathcal{F}^{'}$. 
\end{definition}

\begin{definition}
\label{equ-2}
A \emph{unitary equivalence} between two Fredholm modules over $A$: $(\mathcal{H}, \mathcal{F})$ and $(\mathcal{H}^{'}, \mathcal{F}^{'})$, is given by 
an even unitary isomorphism $T : \mathcal{H} \to \mathcal{H}^{'}$, which intertwines the actions of $A$ and the operators $\mathcal{F}, \mathcal{F}^{'}$. 
\end{definition}

\begin{definition}
The \emph{$K$-homology} $KK_{0}(A,\C)$ is defined as the quotient of the set of Fredholm modules over $A$ by the equivalence relation generated by operator homotopies and unitary equivalences. The addition is given by direct sum and the \emph{zero module} has zero Hilbert space, zero representation, and zero operator.   
\end{definition}

We can also define $KK_{1}(A, \C) = KK_{0}(C_{0}(\R)\otimes A, \C)$ to make $KK_{\bullet}(A, \C)$ into a $\Z_{2}$-graded homology theory of $C^{*}$-algebras. In a special case when $A = \C$, we have that
\[
KK_{0}(\C, \C) \cong \Z, \ \ \ \ \  \ \ KK_{1}(\C, \C) = 0. 
\]
The following two examples are important in this paper. 

\begin{example}
\label{dirac}
Suppose that $(M, g)$ is an even dimensional complete Riemannian manifold. Let $\mathrm{Cliff}(TM)$ be the bundle of complex Clifford algebras over $M$ with standard $\Z_{2}$-grading, and $C_{\tau}(M)$ the $\Z_{2}$-graded $C^{*}$-algebra consisting of all the continuous sections of $\mathrm{Cliff}(TM)$ vanishing at infinity. For any $e \in TM$, let $e^{*} \in T^{*}M$ corresponds to $e$ via $g$. The action of $e$ on the exterior algebra of cotangent bundle $\Lambda^{*}(T^{*}M)$ is given by the formula:
\begin{equation}
\label{clifford mul}
c(e) := \epsilon(e^{*})  - \iota(e),
\end{equation}
where $\epsilon$ is the exterior multiplication and $\iota$ is the contraction. This induces  an action of $C_{\tau}(M)$ on $\Lambda^{*}(T^{*}M)$. 

Let $d$ be the de Rham differential, $d^{*}$ its dual, and $\mathcal{D} = d + d^{*}$. The completeness of $M$ ensures that $\mathcal{D}$ is essentially self-adjoint so that we can apply functional calculus to $\mathcal{D}$. The pair
\begin{equation}
\label{qq-1}
\mathcal{H} = \Gamma_{L^{2}}(M,\Lambda^{*}(T^{*}M)), \hspace{5mm} \mathcal{F}= \frac{\mathcal{D}}{\sqrt{\mathcal{D}^{2}+1}} 
\end{equation}
gives a Fredholm module over $C_{\tau}(M)$ and determines an element
\[
[d_{M}] \in KK_{0}(C_{\tau}(M), \C).
\] 
We call this element the \emph{Dirac element} after Kasparov \cite[Lemma 4.2]{Kasparov88}.
\end{example}

If $M$ is odd dimensional, one can similarly construct the Dirac element which lives in $KK_{1}(C_{\tau}(M),\C)$.

\begin{definition}
\label{spinc}
Let $M$ be a manifold and $V$ an even-dimensional Euclidean vector bundle over $M$. A \emph{Spin$^{c}$-structure} for $V$ is a $\Z_{2}$-graded Hermitian vector bundle $S$ over $M$ together with an even isomorphism of $*$-algebra bundle: 
\[
c: \mathrm{Cliff}(V)  \xlongrightarrow{\cong} \mathrm{End}(S).
\]
In particular, if $V = TM$ has a Spin$^{c}$-structure, then we say $M$ is a \emph{Spin$^{c}$-manifold} and $S$ is its associated \emph{spinor bundle}. 
\end{definition}

\begin{example}
\label{DM}
Let $M$ be an even-dimensional Riemannian manifold, equipped with a  Spin$^{c}$-structure. We denote by  $S_{M}$ its spinor bundle. By choosing a connection $\nabla$ on $S_{M}$, we can define a Spin$^{c}$-Dirac operator
\begin{equation}
\label{dirac operator}
D_{M} = \sum_{i=1}^{\mathrm{dim}(M)}c(e_{i})\nabla_{e_{i}},
\end{equation}
where $\{e_{i}\}$ is an orthonormal basis for $TM$. If $M$ is complete, then we can choose $\nabla$ in a way such that the operator $D_{M}$ is  essentially self-adjoint. By functional calculus, we form
\begin{equation}
\label{qq-2}
\mathcal{H} = \Gamma_{L^{2}}(M,\Lambda^{*}(T^{*}M)), \hspace{5mm} \mathcal{F} = \frac{D_{M}}{\sqrt{D_{M}^{2} + 1}}.
\end{equation} 
Let $C_{0}(M)$ be the $C^{*}$-algebra of continuous functions on $M$ vanishing at infinity, which acts on $\mathcal{H}$ by pointwise multiplication.  The pair $(\mathcal{H}, \mathcal{F})$ in (\ref{qq-2}) gives a Fredholm module over $C_{0}(M)$ and determines an element $[D_{M}] \in KK_{0}(C_{0}(M), \C)$. 

In the case when $M$ is compact, the map $p : M \to \pt$ collapsing $M$ to a point induces a push-forward map:
\[
p_{*} : KK_{\bullet}(C(M), \C) \to KK_{\bullet}(\C, \C).
\]
The integer 
\[
p_{*}([D_{M}]) \in KK_{\bullet}(\C, \C) \cong \Z
\]
equals the index of the Spin$^{c}$-Dirac operator $D_{M}$ \cite[Chapter 11]{Higson00}. 
\end{example}

\begin{remark}
\label{complete}
Suppose that $(M, g)$ is a compact manifold with a regular boundary in the sense that there is a metric preserving diffeomorphism $\phi$ between a neighborhood $V$ of $\partial M$ in $M$ and $\partial M \times [0,1)$. Then $U = \mathrm{Int}(M)$ is a non-complete manifold. However, we can extend the definition of Dirac element to this case. Let $\gamma: U \to \R$ be a smooth function such that for 
\[
\gamma(\phi^{-1}(m, t)) = t^{2}, \hspace{5mm} m \in \partial M, \ t \leq \frac{1}{2}, 
\]
and 
\[
\gamma(\phi^{-1}(m, t))> \frac{1}{4}, \hspace{5mm} m \in \partial M, \ t > \frac{1}{2},
\]
and $\gamma \equiv 1$ on $U \setminus V$. With the rescaled metric $\tilde{g}:= \frac{1}{\gamma^{2}} \cdot g$, $(U, \tilde{g})$ is a complete Riemannian manifold. In addition, the new Clifford action is given by $\tilde{c} = \frac{1}{\gamma} \cdot c$. As the construction in Example \ref{dirac}, we can define the Dirac element 
\[
[d_{U}] \in KK_{\bullet}(C_{\tau}(U), \C).
\]

\end{remark}

The $K$-homology groups $KK_{\bullet}(A, \C)$ are functorial not only  with respect to $C^{*}$-algebra homomorphisms, but more generally with respect to \emph{Morita morphisms} which we shall now describe.

\begin{definition}
If $B$ is a $\Z_{2}$-graded $C^{*}$-algebra, then a (right) \emph{pre-Hilbert $B$-module} is a $\Z_{2}$-graded complex Banach space $\mathcal{E}$ equipped with a (right) $B$-action and a $B$-valued inner product $\langle \  , \  \rangle_{B} : \mathcal{E} \times \mathcal{E} \to B$, which is compatible with the $\Z_{2}$-grading, and linear in the second variable and conjugate-linear in the first variable. Moreover,  the inner product satisfies 
\begin{itemize}
\item
$(\langle \xi, \eta \rangle_{B})^{*} = \langle \eta, \xi \rangle_{B}$;
\item
$\langle \eta, \xi \rangle_{B} b= \langle \eta, \xi  \cdot b \rangle_{B}$;
\item
$\| \xi \|^{2} = \| \langle \xi, \xi \rangle_{B}\|$
\end{itemize}
for all $\xi, \eta \in \mathcal{E}, b \in B$. The norm completion of $\mathcal{E}$ is a \emph{Hilbert $B$-module}.
\end{definition}

The Hilbert $\C$-modules are precisely  Hilbert spaces. 

\begin{definition}
Let $\mathcal{E}$ be a Hilbert $B$-module.  We say that a linear map $T : \mathcal{E} \to \mathcal{E}$ is \emph{adjointable} if there exists a map $T^{*} : \mathcal{E} \to \mathcal{E}$ such that 
\[
\langle T\xi, \eta \rangle_{B} = \langle \xi, T^{*} \eta \rangle_{B}.
\] 
The space of adjointable maps is denoted by $\mathcal{L}_{B}(\mathcal{E})$. 
\end{definition}

Let $A, B$ be two $\Z_{2}$-graded $C^{*}$-algebras. 

\begin{definition}[\cite{Rieffel74}]
A \emph{Hilbert A-B bi-module} is a pair $(\mathcal{E}, \psi)$ in which $\mathcal{E}$ is a Hilbert $B$-module and $\psi : A \to \mathcal{L}_{B}(\mathcal{E})$ is a $*$-representation of $A$ on $\mathcal{E}$. Two Hilbert $A$-$B$ bi-modules $(\mathcal{E}_{i}, \psi_{i}), i = 1, 2$ are called \emph{unitarily equivalent} if there is an isomorphism $U: \mathcal{E}_{1} \to \mathcal{E}_{2}$ preserving the $B$-valued inner product and $\Z_{2}$-grading, and commuting with $\psi$.  In addition, we say that $(\mathcal{E}, \psi)$ is \emph{non-degenerate} if $\psi(A)\mathcal{E} = \mathcal{E}$. \emph{Morita morphisms} between $A$ and $B$ are given by non-degenerate Hilbert $A$-$B$ bi-modules. 
\end{definition}

Suppose that $(\mathcal{E}, \psi)$ is a Hilbert $A$-$B$ bi-module. For any element $[(\mathcal{H}, \mathcal{F})] \in KK_{\bullet}(B,\C)$, where $\mathcal{H}$ is a Hilbert space with left $B$-module structure. Let $\mathcal{E} \otimes_{B} \mathcal{H}$ be the norm completion of the algebraic tensor product of $\mathcal{E}$ and $\mathcal{H}$ over $B$, which is a Hilbert space with a $*$-representation of $A$. It is more elaborate to construct the bounded operator $\hat{\mathcal{F}}$ on $\mathcal{E} \otimes_{B} \mathcal{H}$. We refer to \cite{Kasparov80, Higson00} for the construction. Hence, we define a map 
\begin{equation}
\label{index map}
\phi_{\mathcal{E},*} : KK_{\bullet}(B, \C) \to KK_{\bullet}(A, \C): [(\mathcal{H}, \mathcal{F})] \mapsto [(\mathcal{E} \otimes_{B} \mathcal{H}, \hat{\mathcal{F}})]. 
\end{equation}
In the viewpoint of Kasparov bivariant $K$-theory, $[(\mathcal{E} \otimes_{B} \mathcal{H}, \hat{\mathcal{F}})]$ is the \emph{Kasparov product} of $(\mathcal{H}, \mathcal{F})$ and $(\mathcal{E}, 0)$.

\begin{theorem}\cite{Kasparov80}
$K$-homology is functorial with respect to Morita morphisms. 
\end{theorem}

Let $M$ be an even-dimensional Riemannian Spin$^{c}$-manifold with spinor bundle $S_{M}$. For any $v, w \in \Gamma(S_{M})$, the formula
\[
\langle v, w \rangle := w \cdot v^{*} \in S_{M} \otimes S_{M}^{*} \cong \mathrm{End}(S_{M}), 
\]
defines a $\mathrm{End}(S_{M})$-valued inner product. Moreover, $C_{0}(M)$ acts on $\Gamma(\mathrm{End}(S_{M}))$ by pointwise multiplication. By Definition \ref{spinc}, we have that
\[
\Gamma(\mathrm{End}(S_{M})) \cong C_{\tau}(M), 
\] 
preserving the $\Z_{2}$-grading and involution. Therefore, the spinor bundle $S_{M}$ may be viewed as a Hilbert 
$C_{0}(M)$-$C_{\tau}(M)$ bi-module.

\begin{theorem}[\cite{Kasparov80,Kasparov88}]
\label{index theorem}
If $M$ is a complete Riemannian Spin$^{c}$-manifold with spinor bundle $S_{M}$, then the push-forward map 
\[
\Phi_{S_{M}} : KK_{\bullet}(C_{\tau}(M), \C) \to KK_{\bullet}(C_{0}(M), \C),
\]
maps $[d_{M}]$ to $[D_{M}]$ (see Example \ref{dirac} and \ref{DM}).
\begin{proof}
One can find a proof in \cite[Proposition C.1]{MR2989614}.
\end{proof}
\end{theorem}

\begin{definition}
For a compact Spin$^{c}$-manifold $M$ with spinor bundle $S_{M}$,  we can define the \emph{push-forward index map}:
\begin{equation}
\label{index map-1}
\phi_{S_{M}}:KK_{\bullet}(C_{\tau}(M), \C) \to KK_{\bullet}(\C, \C)
\end{equation} 
as a composition:
\[
KK_{\bullet}(C_{\tau}(M), \C) \xlongrightarrow{\Phi_{S_{M}}} KK_{\bullet}(C(M), \C) \xlongrightarrow{p_{*}} KK_{\bullet}(\C, \C).
\]
\end{definition}

\begin{remark}
Any two spinor bundles $S_{M}, \tilde{S}_{M}$ for $M$ are related by a $\Z_{2}$-graded line bundle 
\[
L = \mathrm{Hom}_{\mathrm{Cliff}(TM)}(S_{M}, \tilde{S}_{M}).
\]
The two maps $\phi_{S_{M}}, \phi_{\tilde{S}_{M}}$  may be different unless the line bundle $L$ can be trivialized. 
\end{remark}

Let $G$ be a compact Lie group. Everything that we have just discussed generalizes to the $G$-equivariant case.
\begin{definition}
A $\Z_{2}$-graded $G$-$C^{*}$-algebra $A$ is a $\Z_{2}$-graded $C^{*}$-algebra equipped with a $*$-homomorphism 
\[
\alpha : G \to \mathrm{Aut}(A): s \mapsto \alpha_{s},
\]
such that $s \mapsto \alpha_{s}$ is strongly continuous and preserving the grading.
\end{definition}

For any $\Z_{2}$-graded $G$-$C^{*}$-algebra $A$, we can define the $G$-equivariant $K$-homology $KK^{G}_{\bullet}(A, \C)$ with straightforward modification. In particular, if $A = \C$ with trivial $G$-action, we have that
\[
KK^{G}_{0}(\C, \C) \cong R(G), \ KK^{G}_{1}(\C, \C) = 0. 
\]
For more details of equivariant $KK$-theory, we refer to \cite[Section 2]{Kasparov88}.

There is a descent operation which transfers the equivariant $KK$-theory to the non-equivariant theory. It involves the notion of crossed products of $C^{*}$-algebras. We will give a rapid review of the basic construction \cite{Williams07}.

If $A$ is a $\Z_{2}$-graded $G$-$C^{*}$-algebra, then the space $C(G, A)$ of continuous $A$-valued functions on $G$ becomes a $\Z_{2}$-graded $*$-algebra with respect to convolution and involution defined by
\[
(f * g)(s) = \int_{G} f(t) \alpha_{t}(g(t^{-1}s))dt,
\]
and
\[
f^{*}(s)  = \alpha_{s}(f(s^{-1}))^{*}. 
\]
A \emph{covariant representation} of $A$ is an even $*$-representation of $A$ on a $\Z_{2}$-graded Hilbert space, together with a compatible unitary representation of $G$ on the same space, that is 
\[
\pi : A \to \mathcal{L}(\mathcal{H}), \hspace{5mm} \pi: G \to U(\mathcal{H}).
\]
It determines a $*$-representation of $C(G, A)$ by the formula:
\[
\langle \pi(f) v, w \rangle = \int_{G}\big\langle \pi(f(s)) \pi(s)v, w \big\rangle ds, \hspace{5mm} v,w \in \mathcal{H}. 
\]
\begin{definition}
We define the \emph{crossed product} $C^{*}(G, A)$ as the completion of $C(G, A)$ under the norm:
\[
\|f\| = \sup\{ \|\pi(f)\| : \pi \ \mathrm{is \ a \ covariant \ representation \ of} \ A \}.
\]

\end{definition}

If $A = \C$, then we get the group $C^{*}$-algebra $C^{*}(G)$. There is a Morita equivalence: 
\[
C^{*}(G) \cong \bigoplus_{\gamma \in \mathrm{Irrep}(G)}\mathrm{End}(V_{\gamma}).
\]
It follows that
\[
KK_{0}(C^{*}(G),\C) \cong \hat{R}(G), \hspace{5mm} KK_{1}(C^{*}(G),\C) =0.
\]

We can extend the construction of crossed product to Morita morphisms between $C^{*}$-algebras. Let $A,B$ be two $\Z_{2}$-graded $G$-$C^{*}$-algebra and $(\mathcal{E}, \psi)$ a Hilbert $A$-$B$ bi-module, equipped with a strongly continuous homomorphism $u : G \to \mathrm{Aut}(\mathcal{E})$ such that
\[
\langle u_{s}(\xi) , u_{s}(\eta) \rangle_{B} = s \cdot (\langle \xi, \eta \rangle_{B}),
\]
and
\[
 u_{s}(\xi \cdot b) = u_{s}(\xi) \cdot (s \cdot b), \hspace{5mm} u_{s}(\psi(a) \xi ) = \psi(s \cdot a) \cdot u_{s}(\xi),
\]
for $\xi, \eta \in \mathcal{E}, a \in A, b \in B, s \in G$. The crossed product $C^{*}(G, \mathcal{E})$ is a Hilbert $C^{*}(G,A)$-$C^{*}(G,B)$ bi-module, defined as the completion of $C(G, \mathcal{E})$ with respect to the $C^{*}(G,B)$-valued inner product:
\[
\langle \xi, \eta \rangle (t) = \int_{G} \langle \xi(s), u_{s} (\eta(s^{-1}t))\rangle ds;
\]
and with left action of $C(G, A)$ on $C(G, \mathcal{E})$ given by
\[
\hat{\psi}(f)(\xi)(t) = \int_{G} \psi(f(s)) u_{s}(\xi(s^{-1}t))ds.
\]
The Hilbert bi-module $C^{*}(G, \mathcal{E})$ gives a Morita morphism between crossed products $C^{*}(G,A)$ and $C^{*}(G,B)$, and induces a map:
\begin{equation}
\hat{\phi}_{\mathcal{E},*} : KK_{\bullet}(C^{*}(G,B), \C) \to KK_{\bullet}(C^{*}(G,A), \C). 
\end{equation}

Let  $M$ be an even-dimensional  Riemannian $G$-manifold. As in Example \ref{dirac}, we consider the pair
\[
\mathcal{H} = \Gamma_{L^{2}}(M,\Lambda^{*}(T^{*}M)), \hspace{5mm} \mathcal{F} = \frac{\mathcal{D}}{\sqrt{\mathcal{D}^{2}+1}}.
\]
For any $f \in C(G, C_{\tau}(M))$, its action on $\mathcal{H}$ is determined by the equation:
\begin{equation}
\label{extra-3}
\langle \pi(f) v, w \rangle = \int_{G}\big\langle c(f(s)) \cdot (s \cdot v), w \big\rangle ds, \hspace{5mm} v,w \in \mathcal{H}, 
\end{equation}
where $c$ is the Clifford action. This extends to an action of $C^{*}(G, C_{\tau}(M))$ on $\mathcal{H}$. One can check that $(\mathcal{H}, \mathcal{F})$ defines a Fredholm module over $C^{*}(G, C_{\tau}(M))$. When there is no confusion, we still denote by 
\[
\big[d_{M}\big] \in KK_{0}(C^{*}(G, C_{\tau}(M)), \C)
\] 
the element determined by $(\mathcal{H}, \mathcal{F})$. Moreover, if $M$ is compact and carries a $G$-equivariant Spin$^{c}$-structure, then we can extend the push-forward index map:
\begin{equation}
\label{formal index map}
\hat{\phi}_{S_{M}}: KK_{\bullet}(C^{*}(G, C_{\tau}(M)), \C) \to KK_{\bullet}(C^{*}(G), \C) \cong \hat{R}(G).
\end{equation}

\section{Localization of Dirac Element}
Following an idea of Witten \cite{Witten82}, we perturb the de Rham complex by a $G$-invariant vector field and obtain a localization formula for the  Dirac element. In addition, we extend the push-forward index map to non-compact manifolds.

Let $G$ be a compact Lie group and $(M, g)$ a Riemannian $G$-manifold. Let $X$ be a vector field generated by $\g$ over $M$. That is, there exists a $G$-equivariant map $\psi_{X} : M \to \g$ such that 
\begin{equation}
\label{condition-1}
X(m) = \frac{d}{dt} \Big|_{t=0} \mathrm{exp}(-t\psi_{X}(m)) \cdot m, \hspace{5mm} \forall m \in M.
\end{equation} 
We assume that the vanishing set $\{ X = 0\} \subset M$ is compact. Define a $G$-invariant function 
\[
h = \|X\|^{2}: M \to \R
\] 
and pick a regular value $\epsilon> 0$. Then 
\[
U = h^{-1}\big([0, \epsilon)\big) \subset M
\] 
is a $G$-invariant  submanifold with a regular boundary.  As in Remark \ref{complete}, we can rescale the metric by a $G$-invariant function $\gamma : U \to \R$ so that $(U, \tilde{g} = \frac{1}{\gamma^{2}}\cdot g)$ is a complete Riemannian $G$-manifold. The Dirac element 
\[
[d_{U}] \in KK^{G}(C^{*}(G, C_{\tau}(U)), \C),
\] 
is represented by 
\begin{equation}
\label{usual pair}
\mathcal{H} = \Gamma_{L^{2}}(U, \wedge^{*}T^{*}U),  \hspace{5mm} \mathcal{F} = \frac{\mathcal{D}}{\sqrt{\mathcal{D}^{2}+1}}.
\end{equation}
For any $f \in C_{\tau}(M)$, it acts on $\mathcal{H}$ as the composition:
\begin{equation}
\label{extra-4}
C_{\tau}(M) \xrightarrow{\mathrm{restriction}} C_{\tau}(U) \xrightarrow{\frac{1}{\gamma}\cdot c} \mathbb{B}(\mathcal{H}).
\end{equation}
By the formula in (\ref{extra-3}), the above extends to an action of $C^{*}(G, C_{\tau}(M))$ on $\mathcal{H}$. However, the pair (\ref{usual pair}) might not be a Fredholm module over $C^{*}(G, C_{\tau}(M))$ in general. Instead we introduce the perturbed Dirac element.  Let us define a 1-form on $U$ by $X^{*}:=\tilde{g}(X, \cdot)$ and perturb the de Rham differential $d$ and its dual $d^{*}$:
\begin{equation} 
\label{eq-1}
d_{X} = d + \epsilon(X^{*}), \hspace{5mm} d^{*}_{X} = d^{*} + \iota(X).
\end{equation}
Set $\mathcal{D}_{X} = d_{X} + d^{*}_{X}$. 

\begin{lemma}
\label{cartan}
One has that 
\[
\mathcal{D}_{X}^{2} = (d + d^{*})^{2} +  (\mathcal{L}_{X} + \mathcal{L}_{X}^{*}) +  \langle X, X \rangle_{\tilde{g}}.
\]
\begin{proof}
It follows directly from the Cartan homotopy formula: 
\[
d \iota_{X}  + \iota_{X}  d = \mathcal{L}_{X}.
\]
\end{proof}
\end{lemma}
On the rescaled Riemannian manifold $(U, \tilde{g})$, we define a pair
\begin{equation}
\label{perturbed pair}
\mathcal{H} = \Gamma_{L^{2}}(U, \wedge^{*}T^{*}U),  \hspace{5mm} \mathcal{F}_{X} = \frac{\mathcal{D}_{X}}{\sqrt{\mathcal{D}_{X}^{2}+1}},
\end{equation}
by replacing $\mathcal{D}$ with $\mathcal{D}_{X}$.

\begin{lemma}
\label{lem-1}
If the vanishing set $\{X = 0\}$ is compact, then the pair (\ref{perturbed pair}) is a Fredholm module over $C^{*}(G,C_{\tau}(M))$ and determines an element 
\[
[d_{U,X}] \in KK_{\bullet}(C^{*}(G, C_{\tau}(M), \C). 
\]
\begin{proof}
One has to prove that for any $h \in C^{*}(G, C_{\tau}(M))$, the operators 
\[
h \cdot (1-\mathcal{F}^{2}_{X}), \hspace{5mm} [h, \mathcal{F}_{X}]
\] 
are compact on $H$.  We approximate in norm any element $h \in C^{*}(G, C_{\tau}(M))$ by finite sums 
\[
\sum_{i}e_{i} \cdot a_{i} \in C^{*}(G, C_{\tau}(M)),
\]
where $e_{i} \in C(G)$ and $a_{i} \in C_{\tau}(M)$. Because the operator $\mathcal{F}_{X}$ is $G$-invariant,  it is enough to prove that the operators
\begin{equation}
\label{fredholm conditions}
a \cdot e \cdot (1 - \mathcal{F}^{2}_{X}), \hspace{5mm}  e\cdot  [a, \mathcal{F}_{X}]
\end{equation}
are compact. 

Let us consider 
\begin{equation}
\label{formula-0}
 a \cdot (1 + \mathrm{Cas}_{\g}+ \mathcal{D}_{X}^{2})^{-1} =  a \cdot (1 + \mathrm{Cas}_{\g} + \mathcal{D}^{2} + \mathcal{L}_{X}+ \mathcal{L}_{X}^{*} +  \langle X, X \rangle_{\tilde{g}})^{-1},
\end{equation}
where $\mathrm{Cas}_{\g}$ is the Casimir operator. It is enough to show that the above operator is compact after restricting to every isotypic component $\mathcal{H}^{V}$ of $\mathcal{H}$, where $V$ is an irreducible $G$-representation. We notice that both $\mathrm{Cas}_{\g}$ and $(\mathcal{L}_{ X} + \mathcal{L}_{X}^{*})$ are bounded on $\mathcal{H}^{V}$,  and the function 
\[
\langle X(m), X(m) \rangle_{\tilde{g}} = \frac{1}{\gamma^{2}} \cdot \langle X(m), X(m) \rangle_{g}
\] 
tends to infinity as $m$ tends to the boundary ($\gamma$ is the rescaling function). By \cite[Section 3]{Gromov83},
\[
 a \cdot (1 + \mathrm{Cas}_{\g}+ \mathcal{D}_{X}^{2})^{-1} 
\]
is a compact operator on $\mathcal{H}^{V}$ and thus it is compact on $\mathcal{H}$ as well.  

On the other hand, we have that
\begin{equation}
\begin{aligned}
&a \cdot (1 +\mathrm{Cas}_{\g} + \mathcal{D}_{X}^{2})^{-1}\cdot e - a \cdot e \cdot (1 +  \mathcal{D}_{X}^{2})^{-1}\\
&=a \cdot (1 + \mathrm{Cas}_{\g} + \mathcal{D}_{X}^{2})^{-1} \cdot ([e, \mathcal{D}_{X}^{2}] +\mathrm{Cas}_{\g}\cdot e) \cdot (1+\mathcal{D}_{X}^{2})^{-1}\\
& =a \cdot  (1 + \mathrm{Cas}_{\g} + \mathcal{D}_{X}^{2})^{-1} \cdot ( \mathrm{Cas}_{\g}\cdot e) \cdot (1+\mathcal{D}_{X}^{2})^{-1}.
\end{aligned}
\end{equation}
In the last equation, $(1+\mathcal{D}_{X}^{2})^{-1}$ is automatically bounded and so is $\mathrm{Cas}_{\g} \cdot e$ (\cite[Lemma 6.5]{Kasparov14}). Therefore, we prove that the operator
\[
a \cdot e \cdot (1 +  \mathcal{D}_{X}^{2})^{-1} = a \cdot e \cdot (1-\mathcal{F}^{2}_{X})
\]
is compact. One can also show that $e\cdot  [a, \mathcal{F}_{X}]$ is compact by the same argument  in \cite[Lemma 4.2]{Kasparov88}.
\end{proof}
\end{lemma}
For the perturbed Dirac element, we have the following:
\begin{theorem}[Localization of Dirac element]
\label{localization theorem}
Let $(M, g)$ be a complete Riemannian $G$-manifold. Suppose that $X$ is a vector field generated by $\g$ over $M$ and $\{X = 0\}$ is compact in $M$. We have that
\[
[d_{M}] = [d_{U, X}] \in KK_{\bullet}(C^{*}(G, C_{\tau}(M)), \C),
\]
where $U$ is a $G$-invariant open neighborhood of $\{ X=0\}$.
\begin{proof}
By definition, the Dirac element $[d_{M}]$ for $(M, g)$ is represented by 
\begin{equation}
\label{element-1}
\mathcal{H}= \Gamma_{L^{2}}(M, \wedge^{*}T^{*}M),  \hspace{5mm} \mathcal{F}  = \frac{\mathcal{D}}{\sqrt{\mathcal{D}^{2}+1}},
\end{equation}
where $\mathcal{D} = d+ d^{*}$. By homotopy, we can also define $[d_{M}]$ using  
\begin{equation}
\label{element-1}
\mathcal{H} = \Gamma_{L^{2}}(M, \wedge^{*}T^{*}M),  \hspace{5mm} \mathcal{F}_{X} = \frac{\mathcal{D}_{X}}{\sqrt{\mathcal{D}_{X}^{2}+1}},
\end{equation}
where $\mathcal{D}_{X} = d_{X}+ d^{*}_{X}$ is the perturbed differential operator with respect to the metric $g$. 

The perturbed Dirac element $[d_{U, X}]$ for $(U, \tilde{g})$ is given by the pair
\begin{equation}
\label{element-2}
\tilde{\mathcal{H}} = \Gamma_{L^{2}}(U, \wedge^{*}T^{*}U),  \hspace{5mm} \tilde{\mathcal{F} }_{X} = \frac{\tilde{\mathcal{D}}_{X}}{\sqrt{\tilde{\mathcal{D}}_{X}^{2}+1}},
\end{equation}
where $\tilde{\mathcal{D}}_{X} = d_{X}+ d^{*}_{X}$ is the perturbed differential operator with respect to the rescaled metric $\tilde{g}$. For any irreducible $G$-representation $V$, we will show that on the  isotypic $V$-component, 
\begin{equation}
\label{eq-20}
\big[(\mathcal{H}^{V}, \mathcal{F}_{X}|_{\mathcal{H}^{V}}) \big]= \big[ (\tilde{\mathcal{H}}^{V}, \tilde{\mathcal{F} }_{X}|_{\tilde{\mathcal{H}}^{V}})\big].
\end{equation}

We choose a relatively compact subset $U^{'}$ such that $\{ X= 0 \} \subseteq U^{'} \subset U$ and
\[
\langle X(m), X(m) \rangle_{\tilde{g}} \gg (\mathcal{L}_{X} + \mathcal{L}^{*}_{X})|_{\mathcal{H}^{V}}, \hspace{5mm} m \in U\setminus U^{'}. 
\]
Without loss of generality, we can also assume that
\[
\langle X(m), X(m) \rangle_{g} \gg (\mathcal{L}_{X} + \mathcal{L}^{*}_{X})|_{\mathcal{H}^{V}}, \hspace{5mm} m \in M \setminus U^{'}. 
\] 
In fact, we can always multiply $X$ by some positive function $f$. Under the above assumptions, the two operators $\mathcal{F}_{X}|_{\mathcal{H}^{V}}$ and $\tilde{\mathcal{F}}|_{\tilde{\mathcal{H}}^{V}}$ are invertible outside $U^{'}$. That is, there exists a positive constant $C$ such that 
\[
\|\mathcal{F}_{X} u\| \geq C \cdot \|u\|, \hspace{5mm} \|\tilde{\mathcal{F}}_{X}v\| \geq C \cdot \|v\|
\]
 for any $u \in \mathcal{H}^{V}$ with $\mathrm{Supp}$ $u \subset M \setminus U^{'}$, and $v \in \tilde{\mathcal{H}}^{V}$ with  $\mathrm{Supp}$ $v \subset U \setminus U^{'}$. On the other hand the two operators $\mathcal{F}_{X}$ and $\tilde{\mathcal{F}}_{X}$ are homotopic over the relatively compact subset $U^{'}$. Hence, we proved (\ref{eq-20}).

\end{proof}
\end{theorem}

Next we generalize the push-forward index to non-compact manifolds. Suppose that  $M$ is a possibly non-compact Spin$^{c}$ $G$-manifold and $X$ is a $G$-invariant vector field defined in (\ref{condition-1}). We assume that the vanishing set $\{ X = 0\}$ is compact in $M$.
Choose a $G$-invariant, relatively compact open subset $U \subset M$, containing $\{ X = 0\}$. We denote by $S_{M}$ the $G$-equivariant spinor bundle for $M$ and $S_{U}$ the one for $U$ obtained by restriction.

\begin{lemma}
There exists a $G$-equivariant embedding $j : U \hookrightarrow N$ into a compact Spin$^{c}$ $G$-manifold $N$, whose $G$-equivariant spinor bundle $S_{N}$ can be identified with $S_{U}$ after restricting to $U$.
\begin{proof}
Let $\rho$ be a $G$-invariant function on $M$ with compact support, such that $0\leq \rho \leq 1$ and $\rho = 1$ on $U$. We define $\psi : M \times \R \to \R: (x, t) \mapsto \rho(x) - t^{2}$. 
Pick a regular value $0 <c < 1$ of $\psi$, which always exists by the Sard's theorem. The set $N = \psi^{-1}(c)$ is a compact $G$-invariant sub-manifold of $U \times \R$. Since $M$ has a $G$-equivariant Spin$^{c}$-structure, so does $M\times \R$ and $N$. The embedding $j : U \hookrightarrow N$ given by $x \mapsto (x, \sqrt{1 - c})$ has the desired properties. 
\end{proof}
\end{lemma}

As in (\ref{perturbed pair}), we consider the perturbed pair:
\begin{equation}
\label{pair-2}
\mathcal{H} = \Gamma_{L^{2}}(U, \wedge^{*}T^{*}U),  \hspace{5mm} \mathcal{F}_{X} = \frac{\mathcal{D}_{X}}{\sqrt{\mathcal{D}_{X}^{2}+1}}.
\end{equation} 
Let $C^{*}(G, C_{\tau}(N))$ acts on $\mathcal{H}$ factoring through the restriction to $U$ as in (\ref{extra-4}). By Lemma \ref{lem-1}, the pair $(\mathcal{H}, \mathcal{F}_{X})$ defines an element in  $KK_{\bullet}(C^{*}(G, C_{\tau}(N)),\C)$. Since $N$ is compact, we define the \emph{push-forward index} of $M$ to be
\begin{equation}
\label{index-1}
\hat{\phi}_{S_{N}}([\mathcal{H}, \mathcal{F}_{X}]) \in KK_{\bullet}(C^{*}(G),\C) \cong \hat{R}(G).
\end{equation}

\begin{definition}
By the excision property and functoriality of $K$-homology, the element in (\ref{index-1}) does not depend on the choices of $U$, $N$ and the embedding $j$ but only the Spin$^{c}$-structure on $M$ and the vector field $X$. We denote it by 
\[
\hat{\phi}_{S_{M}}[d_{M,X}] \in KK_{\bullet}(C^{*}(G),\C).
\]
\end{definition}

\begin{remark}
The $K$-homology index for  transversally elliptic operator was introduced by P. Julg \cite{Julg82}. In \cite{Kasparov14}, Kasparov discussed a $KK$-theoretic approach to transversally elliptic operators in a more general situation. In addition, using the vector field $X$, Braverman \cite{Braverman02} constructed a transversally elliptic operator $D_{M, X}$ on $M$. The push-forward index 
\[
\hat{\phi}_{S_{M}}[d_{M,X}] \in KK_{\bullet}(C^{*}(G),\C) \cong \hat{R}(G)
\]
equals to the $G$-index of $D_{M,X}$. 
\end{remark}

\section{Dirac Element of Vector Bundles}
The goal of this section is to study the push-forward index of non-compact manifolds which have structures of the total spaces of  vector bundles. Throughout this section, we assume that $T$ is a torus. We denote by $\mathfrak{t}$ its Lie algebra and $V_{\lambda}$ the irreducible representation of $T$ labeled by weight $\lambda \in \mathfrak{t}^{*}$.

\begin{definition}
Given any element $\chi \in \hat{R}(T)$ and a vector $\alpha \in \mathfrak{t}$, we say that $\chi$ is \emph{polarized by $\alpha$} if the multiplicity of $V_{\lambda}$ in $\chi$ is zero whenever $\langle \lambda, \alpha \rangle < 0$; and we say it is \emph{strictly polarized by $\alpha$} if the multiplicity of $V_{\lambda}$ in $\chi$ is zero whenever $\langle \lambda, \alpha \rangle \leq 0$. 
\end{definition}
In particular, if $\chi \in \hat{R}(T)$ is strictly polarized by a vector $\alpha$, then the multiplicity of the trivial representation in $\chi$ must be zero. 

Let $V$ be an Euclidean space on which $T$ acts by isometries and preserves only the origin.  We fix a vector $\alpha \in \mathfrak{t}$. The Lie derivative $\mathcal{L}_{\alpha}$ is skew-adjoint and the operator 
\[
J = \frac{ \mathcal{L}_{\alpha}}{|\mathcal{L}_{\alpha}|} \in \mathrm{End}(V)
\] 
defines a $T$-equivariant almost complex structure on $V$. The spinor bundle 
\[
S_{V} = \wedge^{0,*}(T^{*}V)
\] 
defines a $T$-equivariant Spin$^{c}$-structure for $V$. In addition, we denote by $X_{\alpha}$ the vector field induced by $\alpha \in \mathfrak{t}$. 

\begin{lemma}
\label{vector space}
The push-forward index
\[
\hat{\phi}_{S_{V}}([d_{V, X_{\alpha}}]) \in KK_{\bullet}(C^{*}(T), \C) \cong \hat{R}(T)
\]
is polarized by $\alpha$.
\begin{proof}
One can find an explicit formula for the push-forward index in \cite[Chapter 7]{Guillemin98} and \cite[Proposition 5.4]{Paradan01}.
\end{proof}
\end{lemma}

More generally suppose that $F$ is an even dimensional Spin$^{c}$-manifold and $P$ is a principle bundle over $F$ whose structure group is a compact Lie group $H$. Let us also assume that $H$ acts on $V$ in such way that it commutes with the $T$-action and induces a $H$-action on $S_{V}$. In this case,  the vector bundle $P \times_{H} V \to F$ has a $T$-equivariant Spin$^{c}$-structure, with spinor bundle given by $P \times_{H} S_{V}$.

Moreover  the total space $M$ of the vector bundle obtains a Spin$^{c}$-structure from that of $F$ and $V$ as well. In fact, the tangent bundle $TM$ fits into an exact sequence of vector bundles over $M$:
\[
0 \to \nu \to TM \to \pi^{*}TF \to 0,
\]
where $\nu = P \times_{H} TV$ is the vertical tangent bundle, and $\pi : P \to F$ is the projection map. If we choose a splitting of the sequence then we obtain an isomorphism:
\[
TM \cong \nu \oplus \pi^{*}TF. 
\]
We define
\begin{equation}
\label{spinor de}
S_{M} =  (P \times_{H} S_{V}) \otimes  \pi^{*}S_{F},
\end{equation}
and the tangent vector $v \oplus w \in \nu \oplus \pi^{*}TF$
acts as the operator $v \otimes 1 + 1 \otimes w$
on  (\ref{spinor de}). Hence, $S_{M}$ gives a $T$-equivariant Spin$^{c}$-structure for $M$. We denote by $X_{\alpha}$ be the vector field on $M$ induced by $\alpha \in \mathfrak{t}$ as before.

\begin{theorem}
\label{polarization}
If the base $F$ is compact, then the push-forward index
\[
 \hat{\phi}_{S_{M}}([d_{M,X_{\alpha}}]) \in KK_{\bullet}(C^{*}(T),\C) \cong \hat{R}(T)
\]
is polarized by $\alpha$.
\begin{proof}
One can find an explicit formula for the above push-forward index in \cite[Lecture 3]{Atiyah74}). We now give a simple $K$-homology approach. 
  
Let us consider a special case in which the principle bundle $P$ is trivial: $P = F \times H$ (in this case we might as well take $H = \{e\}$). Then $M = F \times V$ and the Dirac element $[d_{M,X_{\alpha}}]$ can be represented by 
\[
 L^{2}(M, \Lambda^{*}T^{*}M) \cong  L^{2}(F, \Lambda^{*}T^{*}F) \otimes  L^{2}(V, \Lambda^{*}T^{*}V), 
\]
together operator 
\[
\mathcal{D}_{M} = \mathcal{D}_{F} \otimes 1 + 1\otimes \mathcal{D}_{V, X_{\alpha}} 
\]
where $\mathcal{D}_{F} = d + d^{*}$ is the differential operator on $F$ and $\mathcal{D}_{V, X_{\alpha}} = d_{X_{\alpha}} + d^{*}_{X_{\alpha}}$ is the perturbed differential operator on $V$. Since $S_{M} = S_{F} \boxtimes S_{V}$, the push-forward map decomposes as well. Hence, we have that 
\[
\hat{\phi}_{S_{M}}([d_{M,X_{\alpha}}]) = \hat{\phi}_{S_{F}}([d_{F}]) \times \hat{\phi}_{S_{V}}([d_{V,X_{\alpha}}]) \in KK_{\bullet}(C^{*}(T),\C).
\]
Remembering that $T$ acts trivially on $F$, one can see that $\hat{\phi}_{S_{M}}([d_{M,X_{\alpha}}])$ is polarized by $\alpha$ by  Lemma \ref{vector space}.

The proof of the general case is similar. We represent the Dirac element $[d_{M,X_{\alpha}}]$ by the Hilbert space
\begin{equation}
\label{fiber model}
[L^{2}(P,  \Lambda^{*}\pi^{*}(T^{*}F)) \otimes L^{2}(V, \Lambda^{*}T^{*}V)]^{H},
\end{equation}
with operator
\[
\tilde{D} \otimes 1 + 1 \otimes \mathcal{D}_{X_{\alpha}} , 
\]
where  $\tilde{D}$ is the ``lifted" differential operator on $P$, acting on sections of $\Lambda^{*}\pi^{*}(T^{*}F)$ and acting as identity in the vertical $H$-direction. This reduces the problem to the first case. 
\end{proof}
\end{theorem}

\begin{remark}
\label{polarized-1}
Let $\mathcal{S}$ be an arbitrary $T$-equivariant spinor bundle for $M$. It is related with $S_{M}$ in (\ref{spinor de}) by a $T$-equivariant $\Z_{2}$-graded line bundle 
\[
L = \mathrm{Hom}_{\mathrm{Cliff}(TM)}(S_{M}, \mathcal{S}). 
\]
Recall that $T$ acts trivially on the zero section $F$. By restricting to $F$, we have that 
\[
L|_{F} = \tilde{L} \otimes \C_{\beta},
\]
where $\tilde{L}$ is a non-equivariant line bundle on $F$ and $\C_{\beta}$ is the trivial line bundle on which $T$ acts with weight $\beta \in \mathfrak{t}^{*}$. If $\langle \beta, \alpha \rangle > 0$, then the push-forward index
\[
\hat{\phi}_{\mathcal{S}}([d_{M, X_{\alpha}}]) \in KK_{\bullet}(C^{*}(T), \C) \cong \hat{R}(T)
\]
 is strictly polarized by $\alpha$. 
\end{remark}

When the base $F$ is non-compact, we need additional structures. Let $K$ be another compact Lie group with Lie algebra $\mathfrak{k}$ acting on $F$. Suppose that $F$ has a $K$-equivariant Spin$^{c}$-structure and $P$ is a $K$-equivariant principle $H$-bundle over $F$. By the above discussion, the product $M = P \times_{H}V$ is a $K\times T$-equivariant Spin$^{c}$-manifold. Let  $\varphi : F \to \mathfrak{k}$ be a $K$-equivariant map and define
\begin{equation}
\label{map}
\psi : M  \to \mathfrak{k} \oplus \mathfrak{t}: (x, v) \mapsto (\varphi(x), \alpha). 
\end{equation}
The map $\psi$ induces a vector field on $M$, denoted by $X$. Under these assumptions, we can use the product model (\ref{fiber model}) to prove the following theorem.

\begin{theorem}
\label{vanishing}
If  the zero set $\{ X = 0\}$ is contained in a compact subset of the zero section of  $M$, then the push-forward index:
\[
 \hat{\phi}_{S_{M}}([d_{M, X}]) \in KK_{\bullet}(C^{*}(K \times T), \C)
\]
is  polarized by $\alpha$. 
\end{theorem}

\section{The Cubic Dirac Operator}
Dirac operators on homogeneous spaces were discussed in a number of references. We will study a specific algebraically defined operator, that is, the cubic Dirac operator after Kostant \cite{Kostant99}. The key property of the cubic Dirac operator is that its square equals to the quadratic Casimir operator up to low order terms. We use \cite{Meinrenken-book} as our primary reference. 

Let $G$ be a compact Lie group and $\g$ its Lie algebra. We equip $\g$ with an ad-invariant inner product and fix an orthonormal basis:
\[
X_{i},\hspace{5mm} i=1, \dots, \mathrm{dim}(\g).
\] 
Let $S_{\g}$ be a  $\Z_{2}$-graded irreducible representation of the Clifford algebra $\mathrm{Cliff}(\g)$, together with \emph{Clifford action}:
\[
c : \mathrm{Cliff}(\g) \hookrightarrow \mathrm{End}(S_{\g}).
\]
If $\g$ is even dimensional, then $\mathrm{Cliff}(\g) \cong \mathrm{End}(S_{\g})$. The \emph{spin representation} $\mathrm{ad}^{\g}: \g \to \mathrm{End}(S_{\g})$ is defined by the formula:
\[
\mathrm{ad}^{\g}(X) := \frac{1}{4}\sum_{i=1}^{\mathrm{dim}(\g)} c([X, X_{i}])c(X_{i}), \hspace{5mm} X \in \g.
\]
  
\begin{definition}
\label{cubic operator}
Let $V$ be a $G$-representation with infinitesimal $\g$-action $\pi$. For any $q \in \R^{+}$, we define an operator on $V \otimes S_{\g}$ by 
\begin{equation}
\begin{aligned}
D_{\g}^{q} &= \sum_{i=1}^{\mathrm{dim}(\g)} \big( \pi(X_{i}) \otimes c(X_{i}) +  q \otimes  \mathrm{ad}^{\g}(X_{i}) \cdot c(X_{i}) \big)\\
& = \sum_{i=1}^{\mathrm{dim}(\g)} c(X_{i}) \cdot \big(\pi(X_{i}) + q \cdot \mathrm{ad}^{\g}(X_{i})\big).
\end{aligned}
\end{equation}
Note that the second form resembles a geometric Dirac operator for the connection on the spinor bundle $S_{\g}$:
\[
\nabla^{q}_{X} = \pi(X_{i}) + q \cdot \mathrm{ad}^{\g}(X_{i}).
\] 
Taking $q = \frac{1}{2}$, $\nabla^{\frac{1}{2}}$ is the Levi-Civita connection and $D_{\g}^{\frac{1}{2}}$ is the usual Riemannian Dirac operator;
while $q = \frac{1}{3}$ gives the \emph{cubic Dirac operator}, denoted by $D_{\g}$. 
\end{definition}

The cubic Dirac operator  $D_{\g}$ gives an element in the quantum Weil algebra $\mathcal{W}(\g) := \mathcal{U}(\g) \otimes \mathrm{Cliff}(\g)$ \cite{Meinrenken00, Meinrenken-book} and has the following interesting property. 

\begin{theorem}[\cite{Kostant99}]
\label{casimir}
The square of the cubic Dirac operator equals
\[
D_{\g}^{2} = 2 \cdot \mathrm{Cas}_{\g} + \frac{1}{12} \mathrm{tr}(\mathrm{Cas}_{\g}).
\]
\end{theorem}

Let $P_{G,+} \subset \Lambda^{*}$ be the set of dominant weights for $G$. The set of finite dimensional irreducible $G$-representations is parametrized by $P_{G,+}$. In particular, for any irreducible representation $V_{\lambda}$ with highest weight $\lambda \in P_{G,+}$, $D^{2}_{\g}$ acts on $V_{\lambda}\otimes S_{\g}$ as the scalar $\|\lambda+ \rho_{G}\|^{2}$, where $\rho_{G}$ is the half sum of positive roots of $G$. 

More generally, Kostant introduced the \emph{relative cubic Dirac operator} for pairs of  Lie groups.  Let $H \subset G$ be a closed subgroup of the equal rank. Using the inner product we write $\mathfrak{g} = \mathfrak{h} \oplus \mathfrak{p}$. This decomposition induces isomorphisms:
\[
\mathrm{Cliff}({\g}) \cong \mathrm{Cliff}({\h}) \otimes \mathrm{Cliff}({\mathfrak{p}}),
\]
and 
\begin{equation}
\label{spinor module}
S_{\g} \cong S_{\h} \otimes S_{\mathfrak{p}},
\end{equation} 
where $S_{\h}, S_{\p}$ are $\Z_{2}$-graded irreducible modules of $\mathrm{Cliff}({\h})$ and $\mathrm{Cliff}({\p})$. Accordingly the operator $\mathrm{ad}(X)$ breaks up as a sum
\[
\mathrm{ad}^{\g}(X) =\mathrm{ad}^{\h}(X) + \mathrm{ad}^{\p}(X), \hspace{5mm} X \in \h.
\]
Let $D_{\h}^{'}$ be the cubic Dirac operator associated to the $\h$-representation on $V\otimes S_{\g} \cong (V\otimes S_{\p}) \otimes S_{\h}$. 

\begin{definition}
The \emph{relative cubic Dirac operator} $D_{\g,\h}$ is defined as the difference $D_{\g,\h} = D_{\g} - D_{\h}^{'}$, which can be written out as
\begin{equation}
\label{relative cubic operator}
D_{\g, \h} =\sum_{i=1}^{\mathrm{dim}(\p)}\big( \pi(X_{i}) \otimes c(X_{i}) + 1 \otimes \frac{1}{3}\cdot  \mathrm{ad}^{\p}(X_{i}) \cdot c(X_{i}) \big),\end{equation}
where $\{ X_{i} \}_{i=1}^{\mathrm{dim}(\p)}$ is an orthonormal basis for $\p$. 
\end{definition}
For the relative cubic Dirac operator, we have the following theorem. 
\begin{theorem}[\cite{Kostant99}]
\label{relative cubic}
Let $V_{\lambda}$ be an irreducible representation of $G$. On each isotypic $H$-summand  $W_{\mu}$ of $V_{\lambda} \otimes S_{\p}$ labeled by highest weight $\mu \in \Lambda_{H,+}^{*}$, the square of the relative Dirac operator acts as
\[
D_{\g,\h}^{2} \big|_{W_{\mu}} = \|\lambda + \rho_{G}\|^{2} - \|\mu + \rho_{H}\|^{2}. 
\]
\end{theorem}

Suppose that the adjoint $H$-action on $\p$ lifts to a $H$-action on $S_{\p}$. In this case, the bundle 
\[
S_{G/H} = G \times_{H} S_{\p}^{*} \to G/H,
\] 
where $S_{\p}^{*}$ is the dual of the spinor module $S_{\p}$,  gives a $G$-equivariant Spin structure on the homogeneous space $G/H$. Given any $H$-space $W$,  we can form a $G$-equivariant vector bundle $E_{W} = G \times_{H} W$ over $G/H$.  The space of $L^{2}$-sections of the tensor bundle
\[
 S_{G/H}\otimes E_{W} = G \times_{H}(S_{\p}^{*} \otimes W) \to G/H
\]
can be identified with 
\[
\mathcal{H} = \big[ L^{2}(G) \otimes S_{\p}^{*}\otimes W \big]^{H}  \cong \bigoplus_{\lambda \in P_{G,+}} V_{\lambda} \otimes \mathrm{Hom}_{H} (V_{\lambda} \otimes S_{\p}, W),
\]
where the isomorphism comes from the Peter-Weyl theorem. 

We now define a Dirac operator on $\mathcal{H}$. Since the relative cubic Dirac operator $D_{\g,\h}$ is $H$-equivariant, it restrict to an operator $D^{V_{\lambda}}_{\g,\h}$ on 
\[
\mathrm{Hom}_{H} (V_{\lambda} \otimes S_{\p}, W).
\] 
Tensoring each operator $D^{V_{\lambda}}_{\g,\h}$ with the identity operator on $V_{\lambda}$, and summing over $V_{\lambda}$, one obtains an  operator $D_{W}$ on $\mathcal{H}$. 

In addition, the Hilbert space  $\mathcal{H}$ carries a $\Z_{2}$-grading coming from that on $S_{\p}$. We write $\mathcal{H} = \mathcal{H}^{+} \oplus \mathcal{H}^{-}$ and the operator $D_{W}$ anti-commutes with the $\Z_{2}$-grading. Since $D_{W}$ is equivariant, the kernel $\mathrm{ker}(D_{W}) \in R(G)$.

\begin{definition}
We define the \emph{Dirac induction}
\begin{equation}
\label{Dirac induction}
\mathcal{I}_{H}^{G}: R(H) \to R(G) 
\end{equation}
to be the map 
\[
W  \mapsto \mathrm{Index}(D_{W}) = \mathrm{ker}(D_{W}) \cap \mathcal{H}^{+} - \mathrm{ker}(D_{W}) \cap \mathcal{H}^{-}.
\]
\end{definition}
We can describe the Dirac induction map in terms of weights.
\begin{theorem}[\cite{Kostant99}]
\label{dirac map}
Let $W_{\lambda}$ be an irreducible representation of $H$ with highest weight $\lambda \in P_{H,+}$. If there exists a $w \in W_{G}$ in the Weyl group (necessarily unique) such that
\[
\mu + \rho_{G} = w \cdot (\lambda + \rho_{H})
\]
for some $\mu \in P_{G,+}$, then the Dirac induction takes $W_{\lambda}$ to the $(-1)^{w} \cdot V_{\mu}$; Otherwise, the Dirac induction takes $W_{\lambda}$ to $0$. 
\end{theorem}
The Dirac induction in (\ref{Dirac induction}) can be extended to a map:
\[
\hat{\mathcal{I}}_{H}^{G}: \hat{R}(H) \to \hat{R}(G).
\]

\begin{remark}
Landweber \cite{Landweber01} and  Posthuma \cite{Posthuma11} generalized the construction of relative cubic Dirac operators to the loop group setting, in which they obtained different generalizations of the Weyl-Kac  formula. 
\end{remark}

Let us briefly compare with \emph{holomorphic induction}.  Suppose that $\mathcal{O} = G/H$ is an orbit of coadjoint $G$-action on $\g^{*}$ and $U$ is a $H$-equivariant Spin$^{c}$-manifold, with spinor bundle $S_{U}$. It is well-known that $\mathcal{O}$ is a complex $G$-manifold with spinor bundle
\begin{equation}
\label{eq-11}
S_{\mathcal{O}} = G \times_{H} (S_{\p}^{*} \otimes \C_{\rho_{G} -\rho_{H}}),
\end{equation}
where  $\C_{\rho_{G} -\rho_{H}}$ is the complex plane on with $H$ acts with weight $\rho_{G} -\rho_{H}$. 

The product $M = G \times_{H} U$ obtains a  $G$-equivariant Spin$^{c}$-structure from that of $U$ and $\mathcal{O}$. In fact, the tangent bundle $TM$ splits into vertical direction and horizontal direction:
\[
TM \cong   (G \times_{H} TU ) \oplus \pi^{*}T \mathcal{O}. 
\]
where $\pi : M \to \mathcal{O}$ is the projection. We form a tensor product 
\begin{equation}
\label{eq-12}
S_{M} = S_{\nu}\hat{\otimes} \pi^{*}S_{\mathcal{O}},
\end{equation}
where $S_{\nu} = G \times_{H} S_{U}$. The tangent vector 
\[
v \oplus w \in (G \times_{H} TU ) \oplus \pi^{*}T \mathcal{O}
\] 
acts as the operator 
\[
v\hat{\otimes} 1 + 1 \hat{\otimes} w
\]
on  $S_{\nu}\hat{\otimes} \pi^{*}S_{\mathcal{O}}$. Hence, $S_{M}$ defines a spinor bundle for $M$. 

We define two Dirac-type operators on the Hilbert space 
\begin{equation}
\label{hilbert space}
\Gamma_{L^{2}}(M, S_{M}) \cong  \bigoplus_{\lambda \in P_{G,+}} V_{\lambda} \otimes \mathrm{Hom}_{H}\big(V_{\lambda} \otimes S_{\p},  \Gamma_{L^{2}}(U, S_{U}) \otimes \C_{\rho_{G} - \rho_{H}} \big).
\end{equation} 
We first construct a Spin$^{c}$-Dirac operator $D_{M}$ by choosing a connection on the spinor bundle $S_{M}$ as in Example \ref{DM}. 

On the other hand, one can also build a Dirac operator on (\ref{hilbert space}) as a sum
\[
\tilde{D}_{M} =  D_{\mathrm{ver}} + D_{\mathrm{hor}}.
\]
Here the vertical part $D_{\mathrm{ver}}$ acts on the factor $\Gamma_{L^{2}}(U, S_{U})$ as the Spin$^{c}$-Dirac operator $D_{U}$; and the horizontal part $D_{\mathrm{hor}}$ is constructed as in the Dirac induction so that it acts as the relative cubic Dirac operator on the factors $V_{\lambda} \otimes S_{\p}$. One can see that $\tilde{D}_{M}$ is a combination of geometrically defined operator and algebraically defined operator. 

The two operators $\tilde{D}_{M}$ and $D_{M}$ are identical in the vertical direction; while in the horizontal direction, one is the relative cubic Dirac operator and the other is the Spin$^{c}$-Dirac operator. As pointed out in Definition \ref{cubic operator}, these two operators are homotopic.  Therefore, the push-forward index $\hat{\phi}_{S_{M}}([d_{M}]) \in \hat{R}(G)$ can be obtained by induction from  $\hat{\phi}_{S_{U}}([d_{U}]) \in \hat{R}(H)$.  By Theorem \ref{dirac map}, we have the following. 
\begin{corollary}
\label{coro-1}
The push-forward index $\hat{\phi}_{S_{M}}([d_{M}]) \in \hat{R}(G)$ has the same multiplicity at the trivial representation as  $\hat{\phi}_{S_{U}}([d_{U}]) \in \hat{R}(H)$. 
\end{corollary}

\section{Quantization of Hamiltonian $G$-Spaces}
The goal of this section is to re-prove Theorem \ref{[Q,R]=0}. Let us consider a Hamiltonian action by a compact connected Lie group $G$ on a possibly non-compact symplectic manifold $(M, \omega)$.  Suppose that $(M, \omega)$ is pre-quantizable, with pre-quantum line bundle $L$ and connection $\nabla^{L}$. We can define the moment map $\mu : M \to \g^{*}$ by the Kostant's formula:
\[
 \langle \mu, \xi \rangle = \frac{\sqrt{-1}}{2\pi}(\nabla^{L}_{X_{\xi}} - \mathcal{L}_{\xi}), \hspace{5mm} \xi \in \g
\]
where $\mathcal{L}_{\xi}$ is the Lie derivative and $X_{\xi}$ is the vector field induced by the infinitesimal action of $\xi$. 

We choose a $G$-equivariant almost complex structure $J$ so that
\[
g(\cdot, \cdot) = \omega(\cdot, J \cdot)
\]
defines a Riemannian metric on $M$. Coupling with the pre-quantum line bundle $L$, the symplectic 2-form $\omega$ determines a spinor bundle 
\[
S_{M} = \wedge^{0,*}(T^{*}M) \otimes L.
\]

Fix an Ad-invariant inner product of $\g$. Let $X$ be the Hamiltonian vector field of the norm square of the moment map $\|\mu\|^{2}$, which is known as the \emph{Kirwan vector field}. By calculation, we see that 
\[
X = -J(d\| \mu\|^{2})^{*} = 2 \sum_{i=1}^{\mathrm{dim}(\g)}\mu(\xi_{i})\cdot X_{\xi_{i}}
\]
where $\{\xi_{i}\}$ is an orthonormal basis of $\g$. From the viewpoint of (\ref{condition-1}),  the vector field $X$ is induced by $2 \cdot \mu$. The Kirwan vector field plays an essential role in the work of Tian-Zhang \cite{Zhang98} and Paradan \cite{Paradan01}.

Let $U$ be a $G$-invariant small open neighborhood of the vanishing set $\{X = 0\}$. In \cite{Kirwan-book}, Kirwan showed that $\{X = 0\}$ can be decomposed:
\begin{equation}
\label{dec-1}
\{ X = 0\}  = \bigsqcup_{\alpha \in \Gamma} G \cdot (M^{\alpha}\cap \mu^{-1}(\alpha)),
\end{equation}
where $\Gamma$ is a discrete subset of positive Weyl chamber $\mathfrak{t}_{+}$. Accordingly $U$ breaks up into different connected components $U_{\alpha}$. In particular, $U_{0}$ contains $\mu^{-1}(0)$. If $\{X = 0\}$ is compact in $M$, then we have the localization formula:
\begin{equation}
\label{break up}
[d_{M}] = [d_{U,X}] = \sum_{\alpha \in \Gamma} [d_{U_{\alpha},X}] \in KK_{0}(C^{*}(G, C_{\tau}(M)),\C).
\end{equation}
Theorem \ref{[Q,R]=0} follows directly from the following two theorems.

\begin{theorem}
\label{symplectic reduction}
If $0$ is a regular value of $\mu$, then  the multiplicity of the trivial representation in 
\[
\hat{\phi}_{S_{M}|_{U_{0}}}([d_{U_{0},X}]) \in KK_{0}(C^{*}(G),\C) \cong \hat{R}(G)
\]
equals to $Q(M_{0}, \omega_{0}) \in \Z$. 
\end{theorem}

\begin{theorem}
\label{syl}
For each $\alpha  \in \Gamma$, if $\alpha \neq 0$, then the element
\[
\hat{\phi}_{S_{M}|_{U_{\alpha}}}([d_{U_{\alpha},X}]) \in KK_{0}(C^{*}(G),\C) \cong \hat{R}(G)
\]
doesn't contribute to the multiplicity of the trivial representation. 
\end{theorem}

\begin{remark}
For Theorem \ref{symplectic reduction} and \ref{syl}, one can also find different proofs in \cite{Paradan01, Zhang98}.
\end{remark}

We focus on Theorem \ref{syl} for the rest of this section. Let us recall the symplectic cross-section theorem \cite{Guillemin-book}:

\begin{theorem}
Let $(M, \omega)$ be a symplectic manifold with a moment map $\mu : M \to \g^{*}$ arising from a Hamiltonian $G$-action. Let $\alpha$ be a point in $\g^{*}$ and $G_{\alpha}$ its isotropy Lie group. There exists a $G_{\alpha}$-invariant subset $R_{\alpha} \subset \g^{*}$ containing $\alpha$ such that the cross-section $V_{\alpha} = \mu^{-1}(R_{\alpha})$ is a $G_{\alpha}$-invariant symplectic sub-manifold of $M$,and the set
\[
U_{\alpha} :=G \cdot V_{\alpha} \cong G \times_{G_{\alpha}} V_{\alpha}
\]
is diffeomorphic to an open subset in $M$. Moreover, if $M$ has a pre-quantum line bundle $L$, then $V_{\alpha}$ has a pre-quantum line bundle $L_{\alpha}$ as well, with the property that $L|_{U_{\alpha}} \cong G \times_{G_{\alpha}} L_{\alpha}$. The moment map associated to the $G_{\alpha}$-action on $V_{\alpha}$ is obtained by restriction. 
\end{theorem}

We apply the symplectic cross-section theorem to our case and conclude that each component $U_{\alpha}$ has the geometric structure:
\begin{equation}
\label{splitting-s}
U_{\alpha} \cong G \times_{G_{\alpha}} V_{\alpha}.
\end{equation}
Coupling with the pre-quantum line bundle $L_{\alpha}$, the symplectic 2-form on $V_{\alpha}$ determines a $G_{\alpha}$-equivariant  Spin$^{c}$-structure  on $V_{\alpha}$ with spinor bundle  
\[
S_{V_{\alpha}} = \wedge^{0,*}(T^{*}V_{\alpha}) \otimes L_{\alpha}. 
\]
Together with the complex structure on $G/G_{\alpha}$, we obtain a $G$-equivariant Spin$^{c}$-structure on $U_{\alpha}$ as in  (\ref{eq-12}), which is equivalent to that obtained by restricting $S_{M}$  to $U_{\alpha}$. By Corollary \ref{coro-1}, the push-forward index 
\[
\hat{\phi}_{S_{M}|_{U_{\alpha}}}([d_{U_{\alpha},X}]) \in KK_{0}(C^{*}(G), \C) \cong \hat{R}(G).
\] 
has the same multiplicity at the trivial representation as
\[
\hat{\phi}_{S_{V_{\alpha}}}([d_{V_{\alpha},X}]) \in KK_{0}(C^{*}(G_{\alpha}), \C) \cong \hat{R}(G_{\alpha}).
\] 
Therefore, it is enough to show that $\hat{\phi}_{S_{V_{\alpha}}}([d_{V_{\alpha},X}])$ does not contribute the trivial representation. 

The smaller Hamiltonian $G_{\alpha}$-manifold $V_{\alpha}$ has a moment map obtained by restriction. We can apply  the localization process again to $V_{\alpha}$. By induction, we reduce the problem to the case when $G_{\alpha} = G$. We can thus assume $G = K \times T_{\alpha}$, where $T_{\alpha}$ the torus generated by $\alpha$. Accordingly the  moment map breaks up as a sum:
\[
\mu =  \mu_{\mathfrak{k}} + \mu_{\alpha} 
\]
where $\mu_{\mathfrak{k}}, \mu_{\alpha}$ are the composition of $\mu$ with the projection from $\g$ to $\mathfrak{k}$ and $\mathrm{Lie}(T_{\alpha})$ respectively. 

Let $F$ be a $K$-invariant open neighborhood of $M^{\alpha} \cap \mu^{-1}(\alpha)$ in $M^{\alpha}$ and $N$ a $G$-invariant tubular neighborhood of $F$ in $M$. In order to apply the Theorem \ref{vanishing}, we need to deform the moment map $\mu$. We first shrink $F$ and $N$ so that there is a small constant $\epsilon$ satisfying:
\begin{enumerate}
\item 
$\|\alpha -  \mu(n) \| < \epsilon$ for all $n \in N$;
\item
$N^{\gamma}$ is contained in the zero section of $N$ for all $\|\gamma -  \alpha\| < \epsilon, \gamma \in \g$. This is possible because $N^{\alpha}$ preserves only the zero section. 
\end{enumerate}
For $0 \leq t \leq 1$, we define
\[
\mu^{t}(x,v) = \mu_{\mathfrak{k}}(x, t \cdot v) + \mu_{\alpha}(x, v): N \to \mathfrak{g},  
\] 
and $X^{t}$ the vector field induced by  $2 \mu^{t}$. The two conditions above ensures that the vanishing set $\{  X^{t} = 0\}$ is contained in a compact subset of the zero section of $N$ for all $0 \leq t \leq 1$. By homotopy, 
\[
[d_{N, X}]= [d_{N, X^{0}}] \in KK_{0}(C^{*}(G, C_{\tau}(M)), \C). 
\]
By identifying $U_{\alpha}$ with $N$, we have the following:
\begin{proposition}
The element
\[
\hat{\phi}_{S_{M}|_{U_{\alpha}}}([d_{U_{\alpha},X^{0}}]) \in KK_{0}(C^{*}(G), \C) \cong \hat{R}(G)
\]
has zero multiplicity at the trivial representation. 
\begin{proof}
We can now apply Theorem \ref{vanishing} to show that 
\[
 \hat{\phi}_{S_{U_{\alpha}}} ([d_{U_{\alpha},X^{0}}])  \in KK_{0}(C^{*}(K \times T_{\alpha}), \C) \cong \hat{R}(K \times T_{\alpha})
 \] 
 is polarized by $\alpha$. Moreover,  the restriction of the pre-quantum line bundle $L$ to the $T_{\alpha}$-fixed zero section $F$ has a $T_{\alpha}$-weight  $\alpha \in \mathfrak{t}$. Thus, it provides a ``positive factor". By Remark \ref{polarized-1},   $\hat{\phi}_{S_{U_{\alpha}}} ([d_{U_{\alpha},X^{0}}])$ is strictly polarized by $\alpha$. 
\end{proof}
\end{proposition}

\bibliographystyle{plain}
\bibliography{mybib}

\end{document}